# Absolutely superficial sequences

BY NGÔ VIÊT TRUNG

*Viên Toán hoc-Viên Khoa hoc Nghiã Đô, Tù Liêm, Hanoi, Vietnam*

## *Introduction*

Let $A$ be a local ring with maximal ideal $\mathfrak{m}$. Let $M$ be a finitely generated module over $A$. Let $a_1, \ldots, a_r$ be a sequence of elements of $m$. Let $q_i$ denote the ideal $(a_1, \ldots, a_i)$, $i = 1, \ldots, r$, and set $\mathfrak{q}_0 = 0_A$ (the zero ideal of $A$), $\mathfrak{q} = \mathfrak{q}_r$.

*Definition.* $a_1, \ldots, a_r$ is called an *absolutely superficial M-sequence* (abbr. a.s. $M$-sequence) if for each $i = 1, \ldots, r$, $a_i$ is an absolutely superficial element of $\mathfrak{q}$ for the module $M_{i-1} := M/\mathfrak{q}_{i-1}M$, i.e. $(\mathfrak{q}^{n+1}M_{i-1}: a_i) \cap \mathfrak{q}M_{i-1} = \mathfrak{q}^n M_{i-1}$ for all $n > 0$ (cf. ((13), definition 2·1)).

This notion was introduced by P. Schenzel (13) in order to study generalized Cohen–Macaulay (resp. Buchsbaum) modules. All results of (13) concerning a.s. sequences depend heavily on the pecularities of generalized Cohen–Macaulay (resp. Buchsbaum) modules. Hence, at first sight, one might think that the notion of a.s. sequences is formal.

In this paper we shall see that a.s. sequences themselves enjoy many interesting properties relative to different topics of the theory of modules. Our main results may be summarized as follows:

(1) There are various characterizations of a.s. sequences. Some of these characterizations are very simple; e.g. $a_1, \ldots, a_r$ is an a.s. $M$-sequence if and only if $\mathfrak{q}_{i-1}M: a_i^2 = \mathfrak{q}_{i-1}M: \mathfrak{q}$ for $i = 1, \ldots, r$ (Section 1).

(2) A.s. sequences are closely related with other specified sequences of (4), (6), (11), (15). A natural consequence of this fact is the characterization of generalized Cohen–Macaulay (resp. Buchsbaum) modules in terms of a.s. sequences (Section 2).

(3) Graded modules associated with an ideal $\mathfrak{q}$ generated by a a.s. $M$-sequence have simple structures; e.g. the Rees module $R_\mathfrak{q}(M)$ is naturally isomorphic to the symmetric module $S_\mathfrak{q}(M)$, where $R_\mathfrak{q}(M)$ and $S_\mathfrak{q}(M)$ are defined like the Rees algebras and the symmetric algebras of the ring theory, c.f (2). That has some applications in the theory of generalized analytic independence developed in (3), (8), (9), and (17) (Section 3).

(4) For each system of parameters $a_1, \ldots, a_r$ of $M$ there exists a polynomial bounding the Hilbert–Samuel function $l(M/\mathfrak{q}^n M)$, $n \geqslant 0$, and this polynomial coincides with the Hilbert–Samuel polynomial of $l(M/\mathfrak{q}^n M)$ if and only if $a_1, \ldots, a_r$ is a a.s. $M$-sequence. As a consequence, one can estimate the Hilbert–Samuel function $l(M/\mathfrak{a}^n M)$ for an arbitrary ideal $\mathfrak{a}$ of $A$ with $l(M/\mathfrak{a}M) < \infty$ (Section 4).

## 1. *Characterizations*

In this section we will establish the main properties of a.s. sequences.

THEOREM 1·1. *The following conditions are equivalent:*

(i) $a_1, \ldots, a_r$ *is an a.s. M-sequence.*

(ii) *For each* $i = 1, \ldots, r$ *there exists an infinite sequence of positive integers n such that* $[(\mathfrak{q}_{i-1}, \mathfrak{q}^{n+1})M : a_i] \cap \mathfrak{q}M = (\mathfrak{q}_{i-1}, \mathfrak{q}^n)M$.

(iii) $(\mathfrak{q}_{i-1}M : a_i) \cap \mathfrak{q}(a_i, \ldots, a_r)^n M = \mathfrak{q}_{i-1}(a_i, \ldots, a_r)^n M$ *for all* $n \geq 0$ *and* $i = 1, \ldots, r$.

(iv) $(\mathfrak{q}_{i-1}M : a_i) \cap \mathfrak{q}M = \mathfrak{q}_{i-1}M$ *for* $i = 1, \ldots, r$.

(v) $\mathfrak{q}_{i-1}M : a_i^2 = \mathfrak{q}_{i-1}M : \mathfrak{q}$ *for* $i = 1, \ldots, r$.

(vi) $\mathfrak{q}_{i-1}M : a_i^m = \mathfrak{q}_{i-1}M : \mathfrak{q}^n$ *for all* $m, n \geq 1$ *and* $i = 1, \ldots, r$.

(vii) $\mathfrak{q}_{i-1}M : a_i = \bigcup_{n=1}^{\infty} \mathfrak{q}_{i-1}M : \mathfrak{q}^n$ *and* $a_i \notin \mathfrak{p}$ *for all* $\mathfrak{p} \in Ass(M/\mathfrak{q}_{i-1}M) \setminus V(\mathfrak{q})$, $(i = 1, \ldots, r,)$ *where* $V(\mathfrak{q})$ *denotes the set of primes of A containing* $\mathfrak{q}$.

*Proof.* We prove the equivalence of the above conditions by the following diagram:

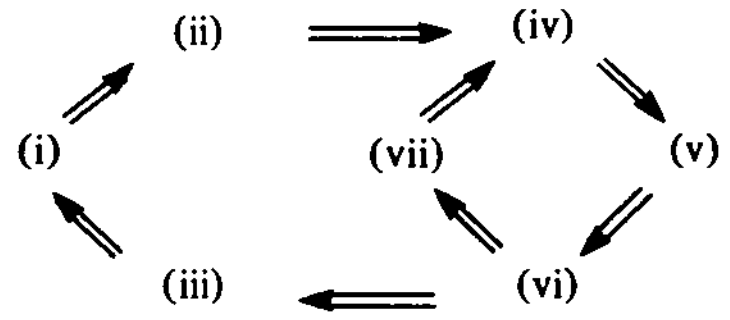

(i) ⇒ (ii) follows from the definition of a.s. sequences.

(ii) ⇒ (iv). We have

$$\mathfrak{q}_{i-1}M \subseteq (\mathfrak{q}_{i-1}M : a_i) \cap \mathfrak{q}M \subseteq \bigcap_n [(\mathfrak{q}_{i-1}, \mathfrak{q}^{n+1})M : a_i] \cap \mathfrak{q}M = \bigcap_n (\mathfrak{q}_{i-1}, \mathfrak{q}^n)M = \mathfrak{q}_{i-1}M,$$

hence (iv).

(iv) ⇒ (v). Dividing both sides of the relation of (iv) by $\mathfrak{q}$ or $a_i$, we get

$$\mathfrak{q}_{i-1}M : a_i\mathfrak{q} = \mathfrak{q}_{i-1}M : \mathfrak{q}, \quad \mathfrak{q}_{i-1}M : a_i^2 = \mathfrak{q}_{i-1}M : a_i,$$

hence (v) because

$$\mathfrak{q}_{i-1}M : q \subseteq \mathfrak{q}_{i-1}M : a_i \subseteq \mathfrak{q}_{i-1}M : a_i\mathfrak{q} \subseteq \mathfrak{q}_{i-1}M : a_i^2.$$

(v) ⇒ (vi). Since

$$\mathfrak{q}_{i-1}M : \mathfrak{q} \subseteq \mathfrak{q}_{i-1}M : a_i \quad (\text{resp. } \mathfrak{q}_{i-1}M : q^2 \subseteq \mathfrak{q}_{i-1}M : a_i^2)$$

we have

$$\mathfrak{q}_{i-1}M : a_i = \mathfrak{q}_{i-1}M : a_i^2 = \mathfrak{q}_{i-1}M : \mathfrak{q} = \mathfrak{q}_{i-1}M : \mathfrak{q}^2,$$

hence (vi).

(vi) ⇒ (vii). It suffices to note that

$$\bigcup_{n=1}^{\infty} \mathfrak{q}_{i-1}M : a_i^n = \bigcup_{n=1}^{\infty} \mathfrak{q}_{i-1}M : \mathfrak{q}^n,$$

iff $a_i \notin \mathfrak{p}$ for all $\mathfrak{p} \in \mathrm{Ass}(M/\mathfrak{q}_{i-1}M) \setminus V(\mathfrak{q})$.

(vii) ⇒ (iv). We have

$$\mathfrak{q}_{i-1}M : a_i = \bigcup_{n=1}^{\infty} \mathfrak{q}_{i-1}M : \mathfrak{q}^n = \bigcup_{n=1}^{\infty} \mathfrak{q}_{i-1}M : a_i^n,$$

hence $\mathfrak{q}_{i-1}M : a_i = \mathfrak{q}_{i-1}M : a_i^2$. From this we can easily deduce that

$$(\mathfrak{q}_{i-1}M : a_i) \cap \mathfrak{q}_iM = \mathfrak{q}_{i-1}M.$$

As a consequence, $(\mathfrak{q}_{r-1}M:a_r)\cap\mathfrak{q}M = \mathfrak{q}_{r-1}M$. If $i<r$, using descending induction on $i$, we have

$$(\mathfrak{q}_{i-1}M:a_i)\cap\mathfrak{q}M = (\mathfrak{q}_{i-1}M:a_i)\cap\left(\bigcup_{n=1}^{\infty}\mathfrak{q}_iM:\mathfrak{q}^n\right)\cap\mathfrak{q}M = (\mathfrak{q}_{i-1}M:a_i)\cap\mathfrak{q}_iM = \mathfrak{q}_{i-1}M.$$

(vi) $\Rightarrow$ (iii). We will first show that $(\mathfrak{q}_{i-1}M:a_i)\cap\mathfrak{q}_jM = (\mathfrak{q}_{i-1}M:a_i)\cap\mathfrak{q}_{j-1}M$ for all $j=i,\dots,r$. Let $u$ be an arbitrary element of $(\mathfrak{q}_{i-1}M:a_i)\cap\mathfrak{q}_jM$. Write $u=v+a_jw$ for some $v\in\mathfrak{q}_{j-1}M$, $w\in M$. Then

$$v+a_jw\in\mathfrak{q}_{i-1}M:a_i = \mathfrak{q}_{i-1}M:\mathfrak{q}\subseteq\mathfrak{q}_{j-1}M:a_j.$$

Thus,

$$w\in\mathfrak{q}_{j-1}M:a_j^2 = \mathfrak{q}_{j-1}M:a_j,$$

hence $u\in q_{j-1}M$, as required. In this way, we have proved $(\mathfrak{q}_{i-1}M:a_i)\cap\mathfrak{q}M = \mathfrak{q}_{i-1}M$; that is (iii) with $n=0$. For $n>0$, let $u$ be an arbitrary element of

$$(\mathfrak{q}_{i-1}M:a_i)\cap\mathfrak{q}(a_i,\dots,a_r)^nM.$$

If $i=r$, write $u=a_r^n(v+a_rw)$ for some $v\in\mathfrak{q}_{r-1}M$, $w\in M$. Then $w\in\mathfrak{q}_{r-1}M:a_r^{n+1} = \mathfrak{q}_{r-1}M:a_r$, hence $u\in\mathfrak{q}_{r-1}a_r^nM$. If $i<r$, write $u=v+a_iw$ for some

$$v\in\mathfrak{q}(a_{i+1},\dots,a_r)^nM,\ w\in\mathfrak{q}(a_i,\dots,a_r)^{n-1}M.$$

Since

$$u\in\mathfrak{q}_{i-1}M:a_i = \mathfrak{q}_{i-1}M:\mathfrak{q}\subseteq\mathfrak{q}_iM:a_{i+1},\quad v\in\mathfrak{q}_iM:a_{i+1}.$$

Hence, by descending induction on $i$, we may assume that $v\in\mathfrak{q}_i(a_{i+1},\dots,a_r)^nM$. Write $v=e+a_if$ for some $e\in\mathfrak{q}_{i-1}(a_{i+1},\dots,a_r)^nM$ and $f\in(a_{i+1},\dots,a_r)^nM$.
Then, since

$$a_i(w+f) = u-e\in\mathfrak{q}_{i-1}M:a_i,\quad w+f\in\mathfrak{q}_{i-1}M:a_i^2 = \mathfrak{q}_{i-1}M:a_i.$$

Thus, by induction on $n$, we may assume that

$$w+f\in\mathfrak{q}_{i-1}(a_i,\dots,a_r)^{n-1}M.$$

From this it follows that

$$u = e+a_i(w+f)\in\mathfrak{q}_{i-1}(a_i,\dots,a_r)^nM,$$

as required.

(iii) $\Rightarrow$ (i). It is easily seen that (iii) may be also formulated for $M/\mathfrak{q}_{i-1}M$, $i=1,\dots,r$. Therefore, we only need to show that $a_1$ is an absolutely superficial element of $\mathfrak{q}$ for $M$, i.e. $(\mathfrak{q}^{n+1}M:a_1)\cap\mathfrak{q}M = \mathfrak{q}^nM$ for all $n\geqslant 1$. If $r=1$, we have

$$(a_1^{n+1}M:a_1)\cap a_1M = (a_1^nM+(O_M:a_1))\cap a_1M = a_1^nM.$$

If $r>1$, we have

$$(a_1M:a_2)\cap(a_2,\dots,a_r)^{n+1}M\subseteq a_1(a_2,\dots,a_r)^nM.$$

Dividing both sides of this inclusion by $a_1$, we get

$$(a_2,\dots,a_r)^{n+1}M:a_1\subseteq(a_2,\dots,a_r)^nM+(O_M:a_1).$$

On the other hand, it is easy to verify that

$$\mathfrak{q}^{n+1}M:a_1 = \mathfrak{q}^nM+((a_2,\dots,a_r)^{n+1}M:a_1).$$

Thus, $(\mathfrak{q}^{n+1}M:a_1)\cap\mathfrak{q}M = (\mathfrak{q}^nM+(O_M:a_1))\cap\mathfrak{q}M = \mathfrak{q}^nM$, as required.

The proof of Theorem 1·1 is complete.

COROLLARY 1·2. *Let $a_1,\dots,a_r$ be a a.s. M-sequence. Then, for each $i = 1,\dots,r$,*
(i) *$a_1,\dots,a_i$ is a a.s. M-sequence.*
(ii) *$\mathfrak{q}_{i-1}M \cap \mathfrak{q}^{n+1}M = \mathfrak{q}_{i-1}\mathfrak{q}^n M$ for all $n \geq 0$.*
(iii) *$(\mathfrak{q}_{i-1},\mathfrak{q}^{n+1})M : a_i = \mathfrak{q}^n M + (\mathfrak{q}_{i-1}M : a_i)$ for all $n \geq 0$.*

*Proof.* (i) (resp. (ii)) follows from the condition (iv) (resp. (iii)) of Theorem 1·1. Since $a_i,\dots,a_r$ is a a.s. $M/\mathfrak{q}_{i-1}M$-sequence for all $i = 1,\dots,r$, (iii) is only a consequence of the fact $\mathfrak{q}^{n+1}M : a_1 = \mathfrak{q}^n M + (O_M : a_1)$, which has been proved in the proof of Theorem 1·1 (iii) ⇒ (i).

*Remark.* The conditions (iv), (v), and (vii) of Theorem 1·1 are very practical in checking whether a given sequence is absolutely superficial or not. For example, let $A = k[[X_1,X_2,X_3]]/(X_1^2, X_1X_2X_3, X_1X_3^2)$, where $k$ is a field. Then, using each of these conditions, one can check that $X_2$, $X_3^2$ is an a.s. $A$-sequence, whereas $X_3^2$, $X_2$ is not. This example shows that the property of being an a.s. sequence is not stable under permutation.

## 2. *Relations to other specified sequences*

In this section we shall see that a.s. sequences are closely related to some specified sequences of the theory of modules.

Let $\mathfrak{a}$ be an arbitrary ideal of $A$.

(1) We call $a_1,\dots,a_r$ an *$\mathfrak{a}$-filter-regular M-sequence* if $a_i \notin \mathfrak{p}$ for all $\mathfrak{p} \in \operatorname{Ass}(M/\mathfrak{q}_{i-1}M)\setminus V(\mathfrak{a})$, $i = 1,\dots,r$. This notion comes from (4), (11) and has led to some interesting results. For instance, $M_\mathfrak{p}$ is a Cohen–Macaulay module and $\dim M_\mathfrak{p} + \dim A/\mathfrak{p} = \dim M$ for all $\mathfrak{p} \in \operatorname{Supp}(M)\setminus\{\mathfrak{m}\}$ if and only if every system of parameters of $M$ is a $\mathfrak{m}$-filter-regular $M$-sequence ((11), Satz 2·5).

By Theorem 1·1 (vii), an a.s. sequence $a_1,\dots,a_r$ is $\mathfrak{a}$-filter-regular if $\sqrt{\mathfrak{q}} \supseteq \mathfrak{a}$. For converse relation, we have the following

PROPOSITION 2·1. *Let $a_1,\dots,a_r$ be an $\mathfrak{a}$-filter-regular M-sequence in $\mathfrak{a}$. Then, for each $n \geq 0$, there exists an ascending sequence of integers $n \leq n_1 \leq \dots \leq n_r$ such that $a_1^{n_1},\dots,a_r^{n_r}$, is a a.s. M-sequence.*

*Proof.* Since $a_i \in \mathfrak{a}$, we can always find an ascending sequence of integers

$$n \leq n_1 \leq \dots \leq n_r$$

such that

$$(a_1^{n_1},\dots,a_{i-1}^{n_{i-1}})M : a_i^{n_i} \supseteq \bigcup_{m=1}^{\infty} (a_1^{n_1},\dots,a_{i-1}^{n_{i-1}})M : \mathfrak{a}^m,$$

$i = 1,\dots,r$. On the other hand, since $a_1,\dots,a_r$ is an $\mathfrak{a}$-filter-regular $M$-sequence iff $a_1,\dots,a_r$ is a regular $M_\mathfrak{p}$-sequence for all $\mathfrak{p} \in \operatorname{Supp}(M/\mathfrak{q}_i M)\setminus V(\mathfrak{a})$ $(i = 1,\dots,r)$, $a_1^{n_1},\dots,a_r^{n_r}$ is also an $\mathfrak{a}$-filter-regular $M$-sequence. Thus, for $i = 1,\dots,r$,

$$\begin{aligned}(a_1^{n_1},\dots,a_{i-1}^{n_{i-1}})M : a_i^{n_i} &\subseteq \bigcup_{m=1}^{\infty} (a_1^{n_1},\dots,a_{i-1}^{n_{i-1}})M : (a_1^{n_1},\dots,a_r^{n_r})^m \\ &= \bigcup_{m=1}^{\infty} (a_1^{n_1},\dots,a_{i-1}^{n_{i-1}})M : \mathfrak{a}^m;\end{aligned}$$

hence

$$(a_1^{n_1},\dots,a_{i-1}^{n_{i-1}})M : a_i^{n_i} = \bigcup_{m=1}^{\infty} (a_1^{n_1},\dots,a_{i-1}^{n_{i-1}})M : (a_1^{n_1},\dots,a_r^{n_r})^m.$$

By Theorem 1·1 (vii), $a_1^{n_1},\dots,a_r^{n_r}$ is an a.s. $M$-sequence.

(2) A superficial element of $\mathfrak{a}$ for $M$ is an element $b \in \mathfrak{a}$ for which there exist integers $c > 0$, $d \geqslant 0$ such that $(\mathfrak{a}^{n+c}M : b) \cap \mathfrak{a}^d M = \mathfrak{a}^n M$ for all sufficiently large $n$. Thus, we call $a_1, \ldots, a_r$ an $\mathfrak{a}$-*superficial M-sequence* if $a_i$ is a superficial element of $\mathfrak{q}$ for $M/\mathfrak{q}_{i-1}M$, $i = 1, \ldots, r$, cf. (6). Superficial elements (hence, by reduction, superficial sequences) have been proved as a useful concept in studying Hilbert–Samuel functions and multiplicities, see, for example, ((18), ch. VIII).

Let $a_1, \ldots, a_r$ be an $\mathfrak{a}$-superficial $M$-sequence. Then we can find integers $c_i > 0$, $d \geqslant 0$ such that

$$[(q_{i-1}, \mathfrak{a}^{n+c_i})M : a_i] \cap \mathfrak{a}^d M = (q_{i-1}, \mathfrak{a}^n)M$$

for all sufficiently large $n$ ($i = 1, \ldots, r$). From this it follows, similarly as in the proof of Theorem 1·1 (ii) ⇒ (iv) ⇒ (v), that

$$(\mathfrak{q}_{i-1}M : a_i) \cap \mathfrak{a}^d M = \mathfrak{q}_{i-1}M$$

and

$$\mathfrak{q}_{i-1}M : a_i = \mathfrak{q}_{i-1}M : a_i\mathfrak{a}^d = q_{i-1}M : \mathfrak{a}^d = \bigcup_{n=1}^{\infty} \mathfrak{q}_{i-1}M : \mathfrak{a}^n.$$

As a consequence, $a_i \notin \mathfrak{p}$ for all $\mathfrak{p} \in \operatorname{Ass}(M/\mathfrak{q}_{i-1}M) \setminus V(\mathfrak{a})$. Hence, $a_1, \ldots, a_r$ is an $\mathfrak{a}$-filter-regular $M$-sequence, cf. ((6), Satz 1·6). In particular, if $d = 0, 1$,

$$\mathfrak{q}_{i-1}M \subseteq (\mathfrak{q}_{i-1}M : a_i) \cap \mathfrak{q}M \subseteq (\mathfrak{q}_{i-1}M : a_i) \cap \mathfrak{a}M = \mathfrak{q}_{i-1}M,$$

hence $a_1, \ldots, a_r$ is an a.s. $M$-sequence by Theorem 1·1 (iv).

(3) Following (11) and (15), we call $a_1, \ldots, a_r$ a $\mathfrak{a}$-*weak M-sequence* if

$$\mathfrak{q}_{i-1}M : a_i \subseteq \mathfrak{q}_{i-1}M : \mathfrak{a} \quad (i = 1, \ldots, r).$$

This notion was used to characterize Buchsbaum (resp. generalized Cohen–Macaulay) modules developed from an answer of W. Vogel to a question of D. A. Buchsbaum. We recall that $M$ is called a Buchsbaum (resp. generaliezd Cohen–Macaulay) module if $l(M/\mathfrak{q}M) - e(\mathfrak{q}; M)$ is (resp. bounded above by) an invariant of $M$ for all parameter ideals $\mathfrak{q}$ of $M$, where $l(M/\mathfrak{q}M)$ denotes the length of $M/\mathfrak{q}M$ and $e(\mathfrak{q}; M)$ is the multiplicity of $M$ relative to $\mathfrak{q}$. This is equivalent to the condition that every system of parameters of $M$ (resp. in $\mathfrak{m}^n$) is a $\mathfrak{m}$-weak $M$-sequence (resp. a $\mathfrak{m}^n$-weak $M$-sequence for some $n \geqslant 1$). See (11), (14), (15), (16) for more informations.

Clearly, every $\mathfrak{a}$-weak sequence is $\mathfrak{a}$-filter-regular. Further, by Theorem 1·1 (vi), an a.s. sequence $a_1, \ldots, a_r$ is $\mathfrak{a}$-weak for all ideals $\mathfrak{a} \subseteq \mathfrak{q}$. For the converse relation, we have the following

Proposition 2·2. *$a_1, \ldots, a_r$ is an a.s. M-sequence if one of the following conditions is satisfied*:

(i) *$a_1, \ldots, a_{i-1}, a_i^2$ is an $\mathfrak{a}$-weak M-sequence ($i = 1, \ldots, r$).*

(ii) *$a_1, \ldots, a_r$ is an $\mathfrak{a}$-weak M-sequence in $\mathfrak{a}^2$.*

*Proof.* Suppose (i). Then

$$\mathfrak{q}_{i-1}M : \mathfrak{q} \subseteq \mathfrak{q}_{i-1}M : a_i^2 \subseteq \mathfrak{q}_{i-1}M : \mathfrak{a} \subseteq \mathfrak{q}_{i-1}M : \mathfrak{q}.$$

Now suppose (ii). Then

$$\mathfrak{q}_{i-1}M : a_i \subseteq \mathfrak{q}_{i-1}M : \mathfrak{a} \subseteq \mathfrak{q}_{i-1}M : \mathfrak{a}^2 \subseteq \mathfrak{q}_{i-1}M : a_i.$$

From this it follows that

$$\mathfrak{q}_{i-1}M:a_i = \bigcup_{m=1}^{\infty} \mathfrak{q}_{i-1}M:\mathfrak{a}^m = \bigcup_{m=1}^{\infty} \mathfrak{q}_{i-1}M:\mathfrak{q}^m.$$

Therefore, the statement follows from Theorem 1·1 (v) and (vii).

PROPOSITION 2·3. *Suppose that $a_1, \ldots, a_r$ is an $\mathfrak{a}$-filter-regular $M$-sequence and there exists a generating set $S$ for $\mathfrak{a}$ such that $a_1, \ldots, a_{r-1}, b$ is an a.s. $M$-sequence for all $b \in S$. Then $a_1, \ldots, a_r$ is an $\mathfrak{a}$-weak $M$-sequence.*

*Proof.* From the first assumption we get

$$\mathfrak{q}_{r-1}M:a_r \subseteq \bigcup_{m=1}^{\infty} \mathfrak{q}_{r-1}M:\mathfrak{a}^m.$$

Using Theorem 1·1 (vi), from the second assumption we get

$$\mathfrak{q}_{r-1}M:b = \bigcup_{m=1}^{\infty} \mathfrak{q}_{r-1}M:b^m \supseteq \bigcup_{m=1}^{\infty} \mathfrak{q}_{r-1}M:\mathfrak{a}^m,$$

and

$$\mathfrak{q}_{i-1}M:a_i = \mathfrak{q}_{i-1}M:(\mathfrak{q}_{r-1}, b) \subseteq \mathfrak{q}_{i-1}M:b$$

for all $b \in S$ ($i = 1, \ldots, r-1$). Thus,

$$\mathfrak{q}_{i-1}M:a_i \subseteq \bigcap_{b \in S} \mathfrak{q}_{i-1}M:b = \mathfrak{q}_{i-1}M:\mathfrak{a}$$

for all $i = 1, \ldots, r$, as required.

Note that if $\dim M = 1$, there always exists a generating set for $\mathfrak{m}^n$, $n \geqslant 1$, whose elements are parameters of $M$, cf. ((7), lemma 3). Then, from Proposition 2·2 and Proposition 2·3 we immediately get the following result of (13), § 2.

COROLLARY 2·4. *$M$ is a Buchsbaum (resp. generalized Cohen–Macaulay) module if and only if every system of parameters of $M$ (resp. in $\mathfrak{m}^n$ for some fixed $n \geqslant 0$) is an a.s. $M$-sequence.*

In particular, the main result of (7) may be reformulated as follows (Buchsbaum modules are hardly characterized by the help of only one system of parameters):

COROLLARY 2·5. *$M$ is a Buchsbaum module if and only if there exists a system of parameters $a_1, \ldots, a_r$ of $M$ in $\mathfrak{m}^2$ and a generating set $S$ for $\mathfrak{m}$ such that $a_1, \ldots, a_i, b_1, \ldots, b_{r-i}$ is an a.s. $M$-sequence for every $r-i$ element subset $b_1, \ldots, b_{r-i}$ of $S$ ($i = 1, \ldots, r$).*

## 3. *Associated graded modules*

In this section we will study graded modules associated with an ideal generated by an a.s. sequence.

Let $G_\mathfrak{q}(M)$ denote the graded module $\bigoplus_{n=0}^{\infty} \mathfrak{q}^n M/\mathfrak{q}^{n+1}M$ over the graded ring $G_\mathfrak{q}(A) = \bigoplus_{n=0}^{\infty} \mathfrak{q}^n/\mathfrak{q}^{n+1}$. It is well-known that $G_\mathfrak{q}(M)$ carries much information on the structure of $M$, see, for example (5). Let $a_1^*, \ldots, a_r^*$ denote the images of $a_1, \ldots, a_r$ in $\mathfrak{q}/\mathfrak{q}^2$ respectively. Let $Q_i$ denote the ideal of $G_\mathfrak{q}(A)$ generated by $a_1^*, \ldots, a_i^*$ ($i = 1, \ldots, r$), and set $Q = Q_r$. Then we have the following result which will make the study of $G_\mathfrak{q}(M)$ easier because it allows the reduction process.

PROPOSITION 3·1. *Let $a_1, \ldots, a_r$ be a a.s. $M$-sequence. Then*

$$G_\mathfrak{q}(M/\mathfrak{q}_i M) \cong G_\mathfrak{q}(M)/Q_i G_\mathfrak{q}(M) \quad \textit{for each} \quad i = 1, \ldots, r.$$

*Proof.* Using Corollary 1·2 (ii), we have

$$G_{\mathfrak{q}}(M/\mathfrak{q}_iM) = \bigoplus_{n=0}^{\infty} (\mathfrak{q}^n, \mathfrak{q}_i)M/(\mathfrak{q}^{n+1}, \mathfrak{q}_i)M \cong \bigoplus_{n=0}^{\infty} q^nM/(q^{n+1}M + q^nM \cap q_iM)$$

$$= \bigoplus_{n=0}^{\infty} q^nM/(q^{n+1}, q_iq^{n-1})M = G_{\mathfrak{q}}(M)/Q_iG_{\mathfrak{q}}(M)$$

(for $n = 0$ we set $\mathfrak{q}_i\mathfrak{q}^{n-1}M = \mathfrak{q}_iM$).

COROLLARY 3·2. *Let $M^*$ denote the localization of $G_{\mathfrak{q}}(M)$ at the maximal graded ideal of $G_{\mathfrak{q}}(A)$. Then $a_1^*, \dots, a_r^*$ is an a.s. $M^*$-sequence if $a_1, \dots, a_r$ is an a.s. $M$-sequence.*

*Proof.* It is easily seen that $(O_{M^*}: a_1^*) \cap QM^* = O_{M^*}$ if $(\mathfrak{q}^{n+1}M: a_1) \cap \mathfrak{q}M = \mathfrak{q}^nM$ for all $n \geqslant 1$. Hence, applying Proposition 3·1, the statement follows from the definition of a.s. sequences.

Other graded modules associated with $M$ and $\mathfrak{q}$ are the Rees module $R_{\mathfrak{q}}(M)$ and the symmetric module $S_{\mathfrak{q}}(M)$. Let $M[X]$ denote the module $M \otimes_A A[X]$ over the polynomial ring $A[X]: = A[X_1, \dots, X_r]$, and consider the elements of $M[X]$ as polynomials over $M$. Then $R_{\mathfrak{q}}(M)$ (resp. $S_{\mathfrak{q}}(M)$) is defined to be the factor module of $M[X]$ by the submodule generated by all (resp. linear) forms of $M[X]$ vanishing at $a_1, \dots, a_r$. Clearly, $R(A_{\mathfrak{q}})$ (resp. $S_{\mathfrak{q}}(A)$) is just the Rees (resp. symmetric) algebra of $\mathfrak{q}$, cf. (2), and $R_{\mathfrak{q}}(M)$ (resp. $S_{\mathfrak{q}}(M)$) may be considered as a graded module over $R_{\mathfrak{q}}(A)$ (resp. $S_{\mathfrak{q}}(A)$). Moreover, we have $G_{\mathfrak{q}}(M) \cong R_{\mathfrak{q}}(M)/\mathfrak{q}R_{\mathfrak{q}}(M)$.

To compute $G_{\mathfrak{q}}(M)$ or $R_{\mathfrak{q}}(M)$ is more difficult than to compute $S_{\mathfrak{q}}(M)$. For this reason one may ask when $R_{\mathfrak{q}}(M) \cong S_{\mathfrak{q}}(M)$. That is the case for example if $a_1, \dots, a_r$ is a regular $M$-sequence, cf. ((2), § 3), ((17), § 1). But below we have a more general result.

THEOREM 3·3. *Let $a_1, \dots, a_r$ be an a.s. $M$-sequence with $\mathfrak{q}_{r-1}M: a_r \neq M$. Then $R_{\mathfrak{q}}(M) \cong S_{\mathfrak{q}}(M)$.*

*Proof.* Let $F \in M[X]$ be an arbitrary form vanishing at $a_1, \dots, a_r$. We have to show that $F = \Sigma F_iG_i$ for some linear forms $F_i \in M[X]$ vanishing at $a_1, \dots, a_r$ and $G_i \in A[X]$. We go by induction on $r$. For $r = 0$ there is nothing to prove. For $r > 0$ we may assume that $t: =$ the degree of $F$ in $X_r > 0$. It suffices to show that the coefficients of all monomials of $F$ containing $X_r^t$ belong to $\mathfrak{q}_{r-1}M: a_r$. For, there exist linear forms $F_i \in M[X]$ vanishing at $a_1, \dots, a_r$ and $G_i \in A[X]$ such that the degree of $F - \Sigma F_iG_i$ in $X_r$ is smaller than $t$.

We shall use a trick from the proof of Lemma 6 (9). Set $I = \bigcup_{n=1}^{\infty} O_A: a_r^n$, $\bar{A} = A/I$ and $\bar{M} = M/IM$. Mark the image of an element or an ideal of $A$ in $\bar{A}$ by an upper line. Then $\bar{a}_r$ is not a zero-divisor of $\bar{A}$. Since by Theorem 1·1 (vii),

$$\mathfrak{q}_{i-1}M: a_i = \bigcup_{n=1}^{\infty} \mathfrak{q}_{i-1}M: \mathfrak{q}^n \subseteq \bigcup_{n=1}^{\infty} \mathfrak{q}_{i-1}M: a_r^n,$$

$$(\mathfrak{q}_{i-1}, I)M: a_i \subseteq \bigcup_{n=1}^{\infty} \mathfrak{q}_{i-1}M: a_r^n \quad (i = 1, \dots, r).$$

Thus,

$$\bar{\mathfrak{q}}_{i-1}\bar{M}: \bar{a}_i \subseteq \bigcup_{n=1}^{\infty} \bar{\mathfrak{q}}_{i-1}\bar{M}: \bar{a}_r^n.$$

Hence $\bar{a}_i \notin P$ for all $P \in \operatorname{Ass}(\bar{M}/\bar{\mathfrak{q}}_{i-1}\bar{M}) \setminus V(\bar{\mathfrak{q}})$. From this it follows that $\bar{a}_1, \ldots, \bar{a}_r$ is a regular $\bar{M}[\bar{a}_r^{-1}]$-sequence, where $\bar{M}[\bar{a}_r^{-1}] := \bar{M} \otimes_{\bar{A}} \bar{A}[\bar{a}_r^{-1}]$. Set $b_1 = \bar{a}_1 \bar{a}_r^{-1}, \ldots, b_{r-1} = \bar{a}_{r-1}\bar{a}_r^{-1}$, $\bar{A}[b] = \bar{A}[b_1, \ldots, b_{r-1}]$, and $\bar{M}[b] = \bar{M} \otimes_{\bar{A}} \bar{A}[b]$. Then, since $\bar{M}[b, \bar{a}_r^{-1}] = M[\bar{a}_r^{-1}]$, $b_1, \ldots, b_{r-1}$ is a regular $\bar{M}[b, \bar{a}_r^{-1}]$-sequence too. On the other hand, we always have

$$\bar{\mathfrak{q}}_{r-1}\bar{M} : \bar{a}_r \subseteq (b_1, \ldots, b_{r-1})\bar{M}[b] \cap \bar{M} \subseteq \bigcup_{n=1}^{\infty} \bar{\mathfrak{q}}_{r-1}\bar{M} : \bar{a}_r^n.$$

But by Theorem 1·1 (vi),

$$\bar{\mathfrak{q}}_{r-1}\bar{M} : \bar{a}_r = \bigcup_{n=1}^{\infty} \bar{\mathfrak{q}}_{r-1}\bar{M} : \bar{a}_r^n.$$

Hence we can conclude that

$$\bar{M}[b]/(b_1, \ldots, b_{r-1})\bar{M}[b] \cong \bar{M}/(b_1, \ldots, b_{r-1})\bar{M}[b] \cap \bar{M} = \bar{M}\Big/\bigcup_{n=1}^{\infty} \bar{\mathfrak{q}}_{r-1}\bar{M} : \bar{a}_r^n,$$

from which it follows that $\bar{a}_r$ is not a zero-divisor on $(b_1, \ldots, b_{r-1})\bar{M}[b]$. So we have proved that the condition $(b)$ of ((5), Satz 3·16·4) is satisfied for the sequence $b_1, \ldots, b_{r-1}$ of $\bar{A}[b]$ and the module $\bar{M}[b]$. Now, write $F = GX_r^t + H$ for some $G \in M[X_1, \ldots, X_{r-1}]$ and $H \in M[X]$ such that the degree of $H$ in $X_r$ is smaller than $t$. Let $\bar{G}$ denote the image of $G$ in $\bar{M}[X_1, \ldots, X_{r-1}]$. Then $\bar{G}(b_1, \ldots, b_{r-1}) \in (b_1, \ldots, b_{r-1})^{s+1}\bar{M}[b]$, where $s$ is the degree of $\bar{G}$. Thus, by the equivalence $(b) \Leftrightarrow (d)$ of (5), Satz 3·16·4, all coefficients of $\bar{G}$ must belong to $(b_1, \ldots, b_{r-1})\bar{M}[b] \cap \bar{M} = \bar{\mathfrak{q}}_{r-1}\bar{M} : \bar{a}_r$. Hence all coefficients of $G$ belong to $(\mathfrak{q}_{r-1}, I)M : a_r = \mathfrak{q}_{r-1}M : a_r$, as required. The proof of Theorem 3·3 is complete.

Theorem 3·3 has some consequences for the theory of generalized analytic independence. Following (3) and (17), we call $a_1, \ldots, a_r$ *N-independent* for $M$, where $N$ is a proper submodule of $M$, if every form of $M[X]$ vanishing at $a_1, \ldots, a_r$ has all its coefficients in $N$.

COROLLARY 3·4. *Let $a_1, \ldots, a_r$ be an a.s. M-sequence. The following conditions are equivalent*:

(i) *$a_1, \ldots, a_r$ are N-independent for M.*

(ii) *$\mathfrak{q}_{r-1}M : a_r \subseteq N$ by every permutation of $a_1, \ldots, a_r$.*

*Moreover, if $a_1, \ldots, a_r$ is a system of parameters of M,* (i) *and* (ii) *are equivalent to*

(iii) $l(\mathfrak{q}^n M/\mathfrak{q}^n N) = \binom{n+r-1}{r-1} l(M/N) < \infty$ *for some (or all)* $n \geq 1$.

*Proof.* By Theorem 3·3, (i)$\Leftrightarrow$(ii) is immediate. Suppose (i) and (ii). Then

$$\bigoplus_{n=0}^{\infty} q^n M/q^n N \cong (M/N)[X]$$

and $N \supseteq \mathfrak{q}M$. Hence (iii) is satisfied if $a_1, \ldots, a_r$ is a system of parameters of $M$. Now suppose (iii). Then every composition series of $\mathfrak{q}^n M/\mathfrak{q}^n N$ has length $\binom{n+r-1}{r-1} l(M/N)$. From this we can conclude that every form of degree $n$ of $M[X]$ vanishing at $a_1, \ldots, a_r$ has all its coefficients in $N$, hence so does every linear form of $M[X]$ vanishing at $a_1, \ldots, a_r$; hence (ii).

In particular, Proposition 2·1 and Corollary 3·4 may be used to construct maximal $N$-independent sets in a given ideal $\mathfrak{a}$ of $A$, cf. (9), § 3. Here we will only demonstrate such a construction by reproving the following characterization of unmixed local rings ((9), proposition 10) which generalizes an answer of (8) to a problem of G. Valla ((17), question 3·9).

THEOREM 3·5. *A is unmixed if and only if dim A is the maximum number of* $\mathfrak{m}^t$-*independent elements in* $\mathfrak{m}^t$ *for all sufficiently large t.*

*Proof.* By (8) (lemma 4 and lemma 5), we only need to show that if $A$ is a complete local ring with $\dim A/\mathfrak{p} = r$ for all $\mathfrak{p} \in \mathrm{Ass}(A)$, then there exist $\mathfrak{m}^t$-independent sets of $r$ elements in $\mathfrak{m}^t$ for all $t \geqslant 1$. Let $a_1, \ldots, a_r$ be a $\mathfrak{m}$-filter-regular $A$-sequence (which always exists). Then, by Proposition 2·1, for each $n \geqslant 0$ there exists an ascending sequence of integers $n \leqslant n_1 \leqslant \ldots \leqslant n_r$ such that $a_1^{n_1}, \ldots, a_r^{n_r}$ is an a.s. $A$-sequence. Note that $a_1, \ldots, a_r$ is also a system of parameters of $A$ and that $n_i \to \infty$ if $n \to \infty$, $i = 1, \ldots, r$. Then, by every permutation of $a_1, \ldots, a_r$

$$\bigcap_{n=1}^{\infty} ((a_1^{n_1}, \ldots, a_{r-1}^{n_{r-1}}) : a_r^{n_r}) \subseteq \bigcap_{n=1}^{\infty} \left( \bigcup_{m=1}^{\infty} (a_1^{n_1}, \ldots, a_{r-1}^{n_{r-1}}) : \mathfrak{m}^m \right)$$

$$= \bigcap_{n=1}^{\infty} \Big( \bigcap_{\mathfrak{p} \in \mathrm{Ass}(A/\mathfrak{q}_{r-1}) \setminus \{\mathfrak{m}\}} (a_1^{n_1}, \ldots, a_{r-1}^{n_{r-1}}) A_{\mathfrak{p}} \cap A \Big) = \bigcap_{\mathfrak{p} \in \mathrm{Ass}(A/\mathfrak{q}_{r-1}) \setminus \{\mathfrak{m}\}} O_{A_\mathfrak{p}} \cap A.$$

It is easily seen that every minimal prime ideal of $A$ is contained in some

$$p \in \mathrm{Ass}(A/q_{r-1}) \setminus \{m\}.$$

From this it follows that

$$\bigcap_{\mathfrak{p} \in \mathrm{Ass}(A/\mathfrak{q}_{r-1}) \setminus \{\mathfrak{m}\}} O_{A_p} \cap A = \bigcap_{\mathfrak{p} \in \mathrm{Ass}(A)} O_{A_p} \cap A = O_A$$

Thus, by ((18), theorem 13, p. 270), there exist $n_1, \ldots, n_r \geqslant t$ such that by every permutation of $a_1, \ldots, a_r$, $(a_1^{n_1}, \ldots, a_{r-1}^{n_{r-1}}) : a_r^{n_r} \subseteq \mathfrak{m}^t$. Hence by Corollary 3·4, $a_1^{n_1}, \ldots, a_r^{n_r}$ form an $\mathfrak{m}^t$-independent set in $\mathfrak{m}^t$.

## 4. *Hilbert–Samuel functions*

In the following we shall see that a.s. sequences of parameters may be also characterized by means of their Hilbert–Samuel functions.

First, we set

$$e_i(q; M) = \begin{cases} l(M/\mathfrak{q}M) - l(\mathfrak{q}_{r-1}M : a_r/(\mathfrak{q}_{r-1}M : a_r) \cap \mathfrak{q}M) & \text{if} \quad i = 0, \\ l(\mathfrak{q}_{r-i}M : a_{r-1+i}/(\mathfrak{q}_{r-i}M : a_{r-i+1}) \cap \mathfrak{q}M) & \\ \qquad - l(\mathfrak{q}_{r-i-1}M : a_{r-i}/(\mathfrak{q}_{r-i-1}M : a_{r-i}) \cap \mathfrak{q}M) & \text{if} \quad 0 < i < r, \\ l(O_M : a_1/(O_M : a_1) \cap \mathfrak{q}M) & \text{if} \quad i = r. \end{cases}$$

THEOREM 4·1. *Let* $a_1, \ldots, a_r$ *be a system of parameters of M. Then*

$$l(M/\mathfrak{q}^{n+1}M) \leqslant \sum_{i=0}^{r} \binom{n+r-i}{r-i} e_i(\mathfrak{q}; M)$$

*for all* $n \geqslant 0$. *Equality holds for an infinite sequence of integers* $n \geqslant 0$ *if and only if* $a_1, \ldots, a_r$ *is an a.s. M-sequence.*

*Proof.* For $r = 0$ there is nothing to prove. For $r > 0$ set $\bar{M} = M/a_1M$. Then we have the exact sequence

$$0 \longrightarrow \mathfrak{q}^{n+1}M : a_1/\mathfrak{q}^nM \longrightarrow M/\mathfrak{q}^nM \xrightarrow{a_1} M/\mathfrak{q}^{n+1}M \longrightarrow \bar{M}/\mathfrak{q}^{n+1}\bar{M} \longrightarrow 0,$$

for all $n \geqslant 0$. From this sequence we get

$$l(\mathfrak{q}^n M/\mathfrak{q}^{n+1}M) = l(\bar{M}/\mathfrak{q}^{n+1}\bar{M}) - l(\mathfrak{q}^{n+1}M : a_1/\mathfrak{q}^n M).$$

Note that $e_i(a_2, \ldots, a_r; \bar{M}) = e_i(q; M)$ if $0 \leqslant i < r-1$, and $e_{r-1}(a_2, \ldots, a_r; \bar{M}) = e_{r-1}(\mathfrak{q}; M) + e_r(\mathfrak{q}; M)$. Then, by induction on $r$, we may assume that

$$l(\bar{M}/\mathfrak{q}^{n+1}\bar{M}) \leqslant \sum_{i=0}^{r-1} \binom{n+r-i-1}{r-i-1} e_i(\mathfrak{q}; M) + e_r(\mathfrak{q}; M).$$

On the other hand, since $(O_M : a_1)/(O_M : a_1) \cap \mathfrak{q}M$ may be considered as a submodule of $(\mathfrak{q}^{n+1}M : a_1)/(\mathfrak{q}^{n+1}M : a_1) \cap \mathfrak{q}M$, we have $l(\mathfrak{q}^{n+1}M : a_1/\mathfrak{q}^n M) \geqslant l((\mathfrak{q}^{n+1}M : a_1) \cap \mathfrak{q}M/\mathfrak{q}^n M) + l(O_M : a_1/(O_M : a_1) \cap \mathfrak{q}M) \geqslant e_r(\mathfrak{q}; M)$ for all $n \geqslant 1$. Thus,

$$l(\mathfrak{q}^n M/\mathfrak{q}^{n+1}M) \leqslant \sum_{i=0}^{r-1} \binom{n+r-i-1}{r-i-1} e_i(\mathfrak{q}; M)$$

for all $n \geqslant 1$. Hence

$$\begin{aligned} l(M/\mathfrak{q}^{n+1}M) &= l(M/\mathfrak{q}M) + \sum_{m=1}^{n} l(\mathfrak{q}^m M/\mathfrak{q}^{m+1}M) \\ &\leqslant \sum_{i=0}^{r} e_i(\mathfrak{q}; M) + \sum_{m=1}^{n} \sum_{i=0}^{r-1} \binom{m+r-i-1}{r-i-1} e_i(\mathfrak{q}; M) \\ &= \sum_{i=0}^{r-1} \sum_{m=0}^{n} \binom{m+r-i-1}{r-i-1} e_i(\mathfrak{q}; M) + e_r(\mathfrak{q}; M) = \sum_{i=0}^{r} \binom{n+r-i}{r-i} e_i(\mathfrak{q}; M). \end{aligned}$$

We have proved the first statement of Theorem 4·1.

Note that if $a_1, \ldots, a_r$ is an a.s. $M$-sequence, we have $(\mathfrak{q}^{n+1}M : a_1) \cap \mathfrak{q}M = \mathfrak{q}^n M$, and, by Corollary 1·2 (iii), $\mathfrak{q}^{n+1}M : a_1 = \mathfrak{q}^n M + (O_M : a_1)$, and hence $l(\mathfrak{q}^{n+1}M : a_1/\mathfrak{q}^n M) = e_r(\mathfrak{q}; M)$ for all $n \geqslant 1$. Then, using induction on $r$ and Theorem 1·1 (ii), we can prove, similarly as above, the second statement of Theorem 4·1.

It is known that if $M$ is a Buchsbaum (resp. generalized Cohen–Macaulay) module, then for every system of parameters $a_1, \ldots, a_r$ of $M$ (resp. in $\mathfrak{m}^n$ for $n$ large enough)

$$l(\mathfrak{q}_{i-1}M : a_i/q_{i-1}M) = \sum_{j=0}^{r-i} \binom{r-i}{j} l(H^j_{\mathfrak{m}}(M)),$$

$i = 1, \ldots, r$, where $H^j_{\mathfrak{m}}(M)$ denotes the $j$th local cohomology module of $M$ with support $\{\mathfrak{m}\}$. Combining this fact with Theorem 1·1 (iv), Corollary 2·4, and Theorem 4·1, we immediately get the following result of (13), § 3.

COROLLARY 4·2. *Let $M$ be a Buchsbaum (resp. generalized Cohen–Macaulay) module. Let $a_1, \ldots, a_r$ be a system of parameters of $M$ (resp. in $m^n$ for $n$ large enough). Then $e_0(\mathfrak{q}; M)$ is the multiplicity $e(q; M)$ of $M$ relative to $q$, $e_i(\mathfrak{q}; M) = \sum_{j=0}^{r-i} \binom{r-i-1}{j-1} l(H^j_{\mathfrak{m}}(M))$, $i = 1, \ldots, r$, where $\binom{r-i-1}{-1} := 0$ if $i \neq r$ and $\binom{-1}{-1} := 1$, and, for all $n \geqslant 0$,*

$$l(M/\mathfrak{q}^{n+1}M) = \binom{n+r}{r} e(\mathfrak{q}; M) + \sum_{i=1}^{r} \sum_{j=0}^{r-i} \binom{n+r-i}{r-i} \binom{r-i-1}{j-1} l(H^j_{\mathfrak{m}}(M)).$$

Theorem 4·1 may be also used to estimate the Hilbert–Samuel function $l(M/\mathfrak{a}^n M)$ of an arbitrary ideal $a$ of $A$ with $l(M/\mathfrak{a}M) < \infty$.

COROLLARY 4·3. *Let $a_1, \ldots, a_r$ be a system of parameters of $M$ in $\mathfrak{a}$. Then, for all $n \geqslant 0$,*

$$l(M/\mathfrak{a}^{n+1}M) \leqslant \sum_{i=0}^{r} \binom{n+r-i-1}{r-i} e_i(\mathfrak{q}; M) + \binom{n+r-1}{r-1} l(M/\mathfrak{a}M + (0_M : \mathfrak{q}^n)).$$

*Equality holds for an infinite sequence of integers $n \geqslant 0$ if and only if the following conditions are satisfied:*

(i) $\mathfrak{q}^n\mathfrak{a}M = \mathfrak{a}^{n+1}M$ *for some* $n \geqslant 1$.

(ii) $a_1, \ldots, a_r$ *is an a.s. $M$-sequence.*

(iii) $\mathfrak{q}_{r-1}M : a_r \subseteq \mathfrak{a}M + \bigcup_{m=1}^{\infty} 0_M : \mathfrak{m}^m$ *by every permutation of* $a_1, \ldots, a_r$.

*Proof.* We have $l(M/\mathfrak{a}^{n+1}M) \leqslant l(M/\mathfrak{q}^n\mathfrak{a}M) = l(M/\mathfrak{q}^nM) + l(\mathfrak{q}^nM/\mathfrak{q}^n\mathfrak{a}M)$. Further, it is easily seen that

$$l(\mathfrak{q}^nM/\mathfrak{q}^n\mathfrak{a}M) \leqslant \binom{n+r-1}{r-1} l(M/\mathfrak{a}M + (0_M : \mathfrak{q}^n)).$$

Hence, using Theorem 4·1 for $l(M/\mathfrak{q}^nM)$, the first statement is immediate. Since $0_M : \mathfrak{q}^n = \bigcup_{m=1}^{\infty} 0_M : \mathfrak{m}^m$ for all sufficiently large $n$, the second statement is only a consequence of the above consideration combined with Theorem 4·1 and Corollary 3·4.

In particular, we have the following improved version of ((10), lemma 1·1) (which generalizes a result of S. Abhyankar on the embedding dimension of a Cohen–Macaulay ring ((1), (1))).

THEOREM 4·4. *Let $M$ be a Buchsbaum module with $\dim M = r$ and multiplicity $e(\mathfrak{m}; M) = e$. Let $s < t$ be arbitrary non-negative integers. Then, for all $n \geqslant 1$,*

$$l(M/\mathfrak{m}^{nt+s}M) \leqslant \binom{n+r-1}{r} t^r e + \sum_{i=1}^{r} \sum_{j=0}^{r-i} \binom{n+r-i-1}{r-i} \binom{r-i-1}{j-1} l(H^j_{\mathfrak{m}}(M))$$
$$+ \binom{n+r-1}{r-1} l(M/\mathfrak{m}^sM + (0_M : \mathfrak{m})).$$

*Moreover, if the residue field $A/\mathfrak{m}$ is infinite, equality holds for some $n \geqslant 1$ if and only if $(a_1^t, \ldots, a_r^t)^n \mathfrak{m}^sM = \mathfrak{m}^{nt+s}M$ for some (or every) system of elements $a_1, \ldots, a_r$ in $\mathfrak{m}\backslash\mathfrak{m}^2$ whose images in $\mathfrak{m}/\mathfrak{m}^2 \subset G_{\mathfrak{m}}(A)$ form a homogeneous system of parameters of $G_{\mathfrak{m}}(M)$.*

*Proof.* Without restriction we may assume that $A/\mathfrak{m}$ is infinite. Then there exist elements $a_1, \ldots, a_r \in \mathfrak{m}\backslash\mathfrak{m}^2$ whose images in $\mathfrak{m}/\mathfrak{m}^2$ form a homogeneous system of parameters of $G_{\mathfrak{m}}(M)$, i.e. $\mathfrak{q}\mathfrak{m}^cM = \mathfrak{m}^{c+1}M$ for some $c \geqslant 0$. Clearly, $a_1, \ldots, a_r$ is also a system of parameters of $M$. By Corollary 4·2, $e_0(\mathfrak{q}; M) = e(\mathfrak{q}; M)$, and it is easily seen that $e(\mathfrak{q}; M) = e$, cf. (12), theorem 1. Thus, applying Corollary 4·2, we have

$$l(M/(a_1^t, \ldots, a_r^t)^n M) = \binom{n+r-1}{r} t^r e + \sum_{i=1}^{r} \sum_{j=0}^{r-i} \binom{n+r-i-1}{r-i} \binom{r-i-1}{j-1} l(H^j_{\mathfrak{m}}(M)).$$

Note that by Corollary 2·4, $a_1^t, \ldots, a_{r-1}^t, a_r$ is an a.s. $M$-sequence for all $t \geqslant 1$. Then, by Theorem 1·1 (vi) and Corollary 1·2 (iii),

$$(a_1^t, \ldots, a_{r-1}^t) M : a_r^t = (a_1^t, \ldots, a_{r-1}^t) M : a_r \subseteq \mathfrak{q}^tM : a_r = \mathfrak{q}^{t-1}M + (0_M : a_r) \subseteq \mathfrak{m}^sM + (0_M : \mathfrak{m}).$$

Hence, by Corollary 3·4,

$$l((a_1^t, \ldots, a_r^t)^n M/(a_1^t, \ldots, a_r^t)^n \mathfrak{q}^sM) = \binom{n+r-1}{r-1} l(M/\mathfrak{m}^sM + (0_M : \mathfrak{m})).$$

Now, since

$$l(M/\mathfrak{m}^{nt+s}M) \leqslant l(M/(a_1^t, \ldots, a_r^t)^n \mathfrak{m}^s M)$$
$$= l(M/(a_1^t, \ldots, a_r^t)^n M) + l((a_1^t, \ldots, a_r^t)^n M/(a_1^t, \ldots, a_r^t)\, \mathfrak{m}^s M),$$

the statements are immediate.

*Remark.* In Theorem 4·4, if every system of parameters of $M$ is a $\mathfrak{m}M$-independent set for $M$ (e.g. $M = A$), we may even assume that $s \leqslant t$ (cf. ((10), lemma 1·1)). In this case, using Theorem 1·1 (iv), it is easily seen that every system of parameters of $M$ is an a.s. $\mathfrak{m}M$-sequence. So, similarly, as in the above proof, we have

$$(a_1^t, \ldots, a_{r-1}^t)\, M : a_r = (a_1^t, \ldots, a_{r-1}^t)\, \mathfrak{m}M : a_r$$
$$\subseteq \mathfrak{q}^t \mathfrak{m}M : a_r \subseteq q^{t-1}\mathfrak{m}M + (O_M : a_r) \subseteq \mathfrak{m}^s M + (O_M : \mathfrak{m}),$$

from which the statements then follow.

*Acknowledgement.* After this paper was sent for publication, I learned that for rings, the notion of absolutely superficial sequences has already been studied under the name of $d$-sequences by C. Huneke. 'On the symmetric and Rees algebra of an ideal generated by a $d$-sequence, *J. Algebra* **62** (1980), 268–275'. To define a $d$-sequence $a_1, \ldots, a_r$ of a ring $A$ he used the condition (v) of Theorem 1·1 (roughly speaking) together with the unimportant condition that $a_1, \ldots, a_r$ form a minimal basis of $(a_1, \ldots, a_r)$. In particular, Theorem 3·3 was already proved with a different method (for the case $M = A$) by Huneke. I would like to thank P. Schenzel and G. Valla for mentioning this.